\documentclass[oneside,english]{amsart}
\usepackage[T1]{fontenc}
\usepackage[latin9]{inputenc}
\usepackage{units}
\usepackage{amsthm}
\usepackage{amssymb}
\usepackage{esint}

\makeatletter
\numberwithin{equation}{section} 
\numberwithin{figure}{section} 
\theoremstyle{plain}
\newtheorem{thm}{Theorem}
  \theoremstyle{plain}
  \newtheorem{lem}[thm]{Lemma}
 \theoremstyle{definition}
 \newtheorem*{defn*}{Definition}
  \theoremstyle{plain}
  \newtheorem{prop}[thm]{Proposition}
  \theoremstyle{remark}
  \newtheorem*{rem*}{Remark}
  \theoremstyle{plain}
  \newtheorem{cor}[thm]{Corollary}

\makeatother

\usepackage{babel}

\begin{document}
\global\long\def\A{{\mathcal{A}}}

\global\long\def\B{{\mathcal{B}}}

\global\long\def\C{{\mathcal{C}}}

\global\long\def\D{\mathcal{D}}

\global\long\def\M{{\mathcal{M}}}

\global\long\def\N{{\mathcal{N}}}

\global\long\def\O{{\mathcal{O}}}

\global\long\def\bbc{{\mathbb{C}}}

\global\long\def\bbp{{\mathbb{P}}}

\global\long\def\bbd{{\mathbb{D}}}

\global\long\def\bbt{{\mathbb{T}}}

\global\long\def\bbr{{\mathbb{R}}}

\global\long\def\bbs{{\mathbb{S}}}

\global\long\def\bbn{{\mathbb{N}}}

\global\long\def\dprime{{\prime\prime}}

\global\long\def\vol#1{\mathrm{vol}(#1)}

\global\long\def\volcn#1{\mathrm{vol_{\bbc^{n}}}(#1)}

\global\long\def\volrn#1#2{\mathrm{vol_{\bbr^{#1}}}(#2)}


\global\long\def\disk#1{#1 \bbd}

\global\long\def\cir#1{#1 \bbt}


\global\long\def\pr#1{\bbp\left( #1 \right)}

\global\long\def\Ex{\mathbb{E}}

\global\long\def\Var#1{\mathrm{Var}(#1)}

\global\long\def\Cov#1{\mathrm{Cov}(#1)}

\global\long\def\ind#1#2{\mathbf{1}_{#1}\left(#2\right)}

\title{Hole Probability for Entire Functions represented by Gaussian Taylor
Series}

\author{Alon Nishry}

\address{School of Mathematical Sciences, Tel Aviv University, Tel Aviv 69978,
Israel }

\email{alonnish@post.tau.ac.il}
\begin{abstract}
Consider the Gaussian entire function\[
f\left(z\right)=\sum_{n=0}^{\infty}\xi_{n}a_{n}z^{n},\]
where $\left\{ \xi_{n}\right\} $ is a sequence of independent and
identically distributed standard complex Gaussians and $\left\{ a_{n}\right\} $
is some sequence of non-negative coefficients, with $a_{0}>0$. We
study the asymptotics (for large values of $r$) of the hole probability
for $f\left(z\right)$, that is the probability $P_{H}\left(r\right)$
that $f\left(z\right)$ has no zeros in the disk $\left\{ \left|z\right|<r\right\} $.
We prove that\[
\log P_{H}\left(r\right)=-S\left(r\right)+o\left(S\left(r\right)\right),\]
where\[
S\left(r\right)=2\cdot\sum_{n\ge0}\log^{+}\left(a_{n}r^{n}\right),\]
as $r$ tends to $\infty$ outside a deterministic exceptional set
of finite logarithmic measure.
\end{abstract}

\thanks{Research supported by the Israel Science Foundation of the Israel
Academy of Sciences and Humanities, grant 171/07.}

\maketitle

\section{Introduction}

\markright{HOLE PROBABILITY FOR ENTIRE FUNCTIONS REP. BY GAUSSIAN TAYLOR SERIES}We
study entire functions of the form\begin{equation}
f(z)=\sum_{n=0}^{\infty}\xi_{n}a_{n}z^{n},\label{eq:f_def}\end{equation}
where $\left\{ \xi_{n}\right\} $ are independent standard complex
Gaussians and $\left\{ a_{n}\right\} $ are (non-negative) deterministic
coefficients, with $a_{0}>0$. It is well known (see \cite{BKPV,Kah})
that, almost surely, the series \eqref{eq:f_def} has an infinite
radius of convergence if and only if the (non-random) Taylor series
$\sum a_{n}z^{n}$ has infinite radius of convergence, that is\[
\lim_{n\to\infty}\frac{\log a_{n}}{n}=-\infty.\]
We use the following notation \begin{eqnarray*}
S(r) & = & 2\cdot\sum_{n\ge0}\log^{+}\left(a_{n}r^{n}\right)\end{eqnarray*}
for the weight of the 'important' coefficients. 

We study the probability of the event where the function $f(z)$ has
no zeros inside $\left\{ |z|<r\right\} $:\[
p_{H}\left(r\right)=-\log\pr{f\left(z\right)\ne0\mbox{ inside }|z|<r},\]
We derive the asymptotics of $p_{H}\left(r\right)$ for large values
of $r$.
\begin{thm}
\label{thm:p_H_raw_asymp}There exists an exceptional set $E\subset\left[1,\infty\right)$
which depends only on the coefficients $a_{n}$ such that $\intop_{E}\frac{dt}{t}<\infty$.
For $r\to\infty$ not belonging to the set $E$\begin{equation}
p_{H}(r)=S(r)+o\left(S(r)\right).\label{eq:p_H_asymp}\end{equation}

\end{thm}
In the case where the coefficients $a_{n}$ satisfy some regularity
conditions, this result was proved in \cite{Ni2}. The general case
obtained here uses ideas that were introduced in \cite{Ni2,ST}, as
well as some new ingredients, including accurate lower bounds for
determinants of some large covariance matrices, which might be of
independent interest. Some exceptional set $E$ in Theorem \ref{thm:p_H_raw_asymp}
seems to be unavoidable if we do not require additional regularity
conditions from the coefficients $a_{n}$. In many cases, it can be
dropped and the error term in \eqref{eq:p_H_asymp} can be shown to
be of order $\O\left(\sqrt{S\left(r\right)}\log S\left(r\right)\right)$,
see Section \ref{sub:proof_strategy} below.

\textit{Acknowledgments:} I would like to thank Fedor Nazarov for
insisting that the restrictions from \cite{Ni2} are unnecessary and
suggesting a way to remove them and Mikhail Sodin for several ideas
and remarks.

\section{Preliminaries\label{sec:Main-Thm-Prelim}}

\subsection{Notation and assumptions}

We denote by $\disk r$ the disk $\left\{ z\,\mid\,|z|\le r\right\} $
and by $\cir r$ its boundary $\left\{ z\,\mid\,|z|=r\right\} $,
with $r\ge1$. The letters $c$ and $C$ denote positive absolute
constants (which can change across lines). We also use the standard
notation:\[
M\left(r\right)=\max_{z\in\disk r}\left|f\left(z\right)\right|.\]
In order to simplify some of the expressions in the paper, we will
assume from now on that \[
a_{0}=1.\]

\subsection{A growth lemma by Hayman}

A set $E\subset\bbr^{+}$ is of finite \textsl{logarithmic measure
}(FILM) if\[
\intop_{E}\,\frac{dt}{t}<\infty.\]
The following lemma is a general lemma about the growth of functions
(taken from \cite{Hay}).
\begin{lem}
\label{lem:H_Growth_Lemma}Suppose that $N(r)$ is a positive increasing
function of $r$ for $r\ge r_{0}.$ Then if $\alpha>0$, and $|h|<N(r)^{-\alpha}$,
we have\[
\left|N\left(re^{h}\right)-N(r)\right|<\alpha N(r),\]
for all $r$ outside a set of FILM.
\end{lem}

\subsection{Further notations, definition of the exceptional set}

We use the following notations:\[
b_{n}=\begin{cases}
\frac{\log a_{n}}{n}+\log r & ,\, a_{n}>0,\\
-\infty & ,\, a_{n}=0,\end{cases}\]
for the regularized coefficients, and \begin{eqnarray*}
N\left(r\right) & = & \left\{ n\,\mid\, b_{n}\ge0\right\} ,\\
S(r) & = & 2\cdot\sum_{n\in N\left(r\right)}\log\left(a_{n}r^{n}\right)=2\cdot{\displaystyle \sum_{n\in N\left(r\right)}}nb_{n}.\end{eqnarray*}
Since the coefficients $a_{n}$ satisfy $\frac{\log a_{n}}{n}\to-\infty$,
we have that $b_{n}\to-\infty$ as $n\to\infty$ and the set $N(r)$
is finite for every $r\in\bbr$.

We now write \begin{eqnarray*}
n\left(r\right) & = & \#N(r),\\
m\left(r\right) & = & \sum_{n\in N\left(r\right)}n,\\
N_{\delta}\left(r\right) & = & \left\{ n\,\mid\, b_{n}\left(r\right)\ge-\delta\right\} ,\\
n_{\delta}\left(r\right) & = & \#N_{\delta}\left(r\right).\end{eqnarray*}
Note that $b_{n}\left(r\right)$ is increasing with $r$ and therefore
$n\left(r\right)$ and $m\left(r\right)$ are increasing functions
of $r$. Also notice the following relations: \begin{eqnarray}
N_{-\delta}\left(r\right) & \subset & N_{0}\left(r\right)=N\left(r\right),\label{eq:N_delta_to_N_r}\\
N_{-\delta}\left(r\right) & = & N\left(re^{-\delta}\right).\nonumber \end{eqnarray}
Applying Lemma \ref{lem:H_Growth_Lemma} to $m\left(r\right)$, with
\[
\delta=m^{-\nicefrac{1}{4}}\left(r\right),\]
we have that outside an exceptional set of FILM:\begin{eqnarray}
m\left(re^{-\delta}\right) & > & \frac{3}{4}\cdot m\left(r\right),\label{eq:m(r)_normal}\\
m\left(re^{\delta}\right) & < & \frac{5}{4}\cdot m\left(r\right).\nonumber \end{eqnarray}
We will use the term \textit{normal} for non-exceptional values of
$r$ (the choice of $\delta$ remains the same throughout the paper).
In general if the inequalities \eqref{eq:m(r)_normal} hold for all
large values of $r$, then there is no exceptional set in Theorem
\ref{thm:p_H_raw_asymp}.

\subsection{Estimates for $S(r)$}

First we find relations between $S\left(r\right)$ and the functions
$m\left(r\right),\, n\left(r\right)$ which appear in the error terms
for the hole probability.
\begin{lem}
\label{lem:S(r)_low_bnd}For normal values of $r$ we have\[
S(r)\ge\frac{3}{2}\cdot\left(m\left(r\right)\right)^{3/4}\ge c\cdot n\left(r\right)^{3/2}.\]
\end{lem}
\begin{proof}
First notice that\begin{equation}
m\left(r\right)=\sum_{n\in N\left(r\right)}n\ge c\cdot\left(n\left(r\right)\right)^{2},\label{eq:m(r)_low_bnd}\end{equation}
since $m\left(r\right)$ is minimal when $N\left(r\right)=\left\{ 0,1,\ldots,n\left(r\right)-1\right\} $.
Now, using the relation \eqref{eq:N_delta_to_N_r} between $N\left(r\right)$
and $N_{-\delta}\left(r\right)$, we get:\begin{eqnarray*}
\frac{S\left(r\right)}{2} & = & \sum_{n\in N\left(r\right)}nb_{n}\left(r\right)\ge\sum_{n\in N_{-\delta}\left(r\right)}nb_{n}\left(r\right)\ge\\
 & \ge & \sum_{n\in N\left(re^{-\delta}\right)}n\delta\ge\delta\cdot m\left(re^{-\delta}\right)\ge\\
 & \ge & \frac{3}{4}\cdot\frac{m\left(r\right)}{\left(m\left(r\right)\right)^{1/4}}\ge c\cdot n\left(r\right)^{3/2}.\end{eqnarray*}

\end{proof}
\newpage{}We now estimate the rate of growth of the function $S\left(r\right)$.
\begin{lem}
\label{lem:S(r)_growth}For $\gamma\in\left(0,\frac{1}{2}\right)$
we have,\[
S\left(\left(1-\gamma\right)r\right)\ge S\left(r\right)-4\gamma\cdot m\left(r\right).\]

\begin{proof}
Write $r^{\prime}=\left(1-\gamma\right)r$ and notice that for $\gamma<\frac{1}{2}$
we have\[
\log\left(1-\gamma\right)\ge-\gamma-\gamma^{2}\ge-2\gamma.\]
It follows that (since $N(r^{\prime})\subset N(r)$)\begin{eqnarray*}
\frac{S(r)-S(r^{\prime})}{2} & = & \sum_{n\in N(r)\backslash N(r^{\prime})}\log\left(a_{n}r^{n}\right)+\sum_{n\in N(r^{\prime})}\left\{ \log\left(a_{n}r^{n}\right)-\log\left[a_{n}\left(r^{\prime}\right)^{n}\right]\right\} \\
 & = & {\displaystyle \sum_{n\in N(r)\backslash N(r^{\prime})}}\log\left(a_{n}r^{n}\right)+\sum_{n\in N\left(r^{\prime}\right)}n\cdot\log\left(\frac{1}{1-\gamma}\right)\\
 & = & {\textstyle \sum_{1}+\sum_{2}}.\end{eqnarray*}
For the first sum we notice that if $n\in N(r)\backslash N(r^{\prime})$
then\begin{eqnarray*}
0 & \ge & \log\left[a_{n}\left(r^{\prime}\right)^{n}\right]\ge\log\left(a_{n}r^{n}\right)-2n\gamma\\
 & \Downarrow\\
\log\left(a_{n}r^{n}\right) & \le & 2n\gamma\end{eqnarray*}
and so\[
{\textstyle \sum_{1}}\le2\gamma\cdot\sum_{n\in N(r)\backslash N(r^{\prime})}n\le2\gamma\cdot\left[m\left(r\right)-m\left(r^{\prime}\right)\right]\]
For the second sum we have\[
{\textstyle \sum_{2}}\le2\gamma\cdot m\left(r^{\prime}\right).\]
Overall we have\[
S\left(r\right)-S\left(r^{\prime}\right)\le4\gamma\cdot m\left(r\right).\]

\end{proof}
\end{lem}

\subsection{Gaussian Distributions\label{sub:Gauss-props}}

Many times we use the fact that if $a$ has a $N_{\bbc}(0,1)$ distribution,
then \begin{equation}
\pr{|a|\ge\lambda}=\exp(-\lambda^{2}),\label{eq:Gaus_prob_large}\end{equation}
and for $\lambda\le1,$

\begin{equation}
\pr{|a|\le\lambda}\in\left[\frac{\lambda^{2}}{2},\lambda^{2}\right].\label{eq:Gaus_prob_small}\end{equation}

\subsection{Strategy of the proof\label{sub:proof_strategy}}

We had the following
\begin{defn*}
The value $r\ge1$ is \textit{normal} if\[
cm\left(r\right)\le m\left(re^{h}\right)\le Cm\left(r\right),\quad\left|h\right|\le\left(m\left(r\right)\right)^{-\nicefrac{1}{4}}.\]
We show that for normal values of $r$\[
p_{H}\left(r\right)\le S\left(r\right)+C\sqrt{m\left(r\right)}\log m\left(r\right)\]
(Proposition \ref{prp:p_H_upp_bnd}) and also that\[
p_{H}\left(r\right)\ge S\left(r\right)-Cn\left(r\right)\log S\left(r\right)\]
(Proposition \ref{prp:p_H_low_bnd}). Combining these bounds with
Lemma \eqref{lem:S(r)_low_bnd}, we get Theorem \ref{thm:p_H_raw_asymp}
(with an error term $\O\left(S\left(r\right)^{\nicefrac{2}{3}+\epsilon}\right)$).
If the coefficients $a_{n}$ are such that each value $r\ge1$ is
normal, then the exceptional set $E$ in Theorem \ref{thm:p_H_raw_asymp}
is not needed. In some cases, where the coefficients $a_{n}$ have
a regular asymptotic behaviour (for instance%
\footnote{In that case we get $p_{H}\left(r\right)=\frac{e^{2}r^{4}}{4}+\O\left(r^{2}\log r\right)$.
Due to a computational error, in \cite{Ni1} we gave the main term
in this asymptotics with the extra factor $3$.%
} $a_{n}=\frac{1}{\sqrt{n!}}$, or more generally, $a_{n}\sim\frac{1}{\Gamma\left(\alpha n+1\right)}$,
with some $\alpha>0$), then $m\left(r\right)\le CS\left(r\right)$
and we obtain\[
p_{H}\left(r\right)=S\left(r\right)+\O\left(\sqrt{S\left(r\right)}\log S\left(r\right)\right),\quad r\to\infty.\]

\end{defn*}

\section{Upper Bound for $p_{H}(r)$ \label{sec:Main-Thm-Upper-Bound}}

In this section we prove the following
\begin{prop}
\label{prp:p_H_upp_bnd}For normal values of $r$, we have \[
p_{H}(r)\le S(r)+C\cdot\sqrt{m\left(r\right)}\log m(r),\]
with $C$ some positive absolute constant.\end{prop}
\begin{rem*}
We note that $r$ is assumed to be large.
\end{rem*}
The simplest case where $f(z)$ has no zeros inside $\disk r$ is
when the constant term dominates all the others. We therefore study
the event $\Omega_{r}$, which is the intersection of the events $\mbox{{\rm \mbox{(i)}}}$,$\mbox{(ii)}$
and $\mbox{(iii)}$, ($C$ will be selected in an appropriate way)\[
\begin{array}{ll}
\mbox{{\rm \mbox{(i)}}}: & |\xi_{0}|\ge C\left(m\left(r\right)\right)^{1/4},\\
\mbox{{\rm \mbox{(ii)}}}: & {\displaystyle \bigcap_{N\left(r\right)\backslash\left\{ 0\right\} }}\mbox{{\rm \mbox{(ii)}}}_{n},\\
\mbox{{\rm \mbox{(iii)}}}: & {\displaystyle \bigcap_{n\in\widetilde{N}_{\delta}\left(r\right)\backslash N\left(r\right)}}\mbox{{\rm \mbox{(iii)}}}_{n},\\
\mbox{{\rm \mbox{(iv)}}}: & {\displaystyle \bigcap_{n\in\left(\widetilde{N}_{\delta}\left(r\right)\right)^{c}}}\mbox{{\rm \mbox{(iv)}}}_{n},\end{array}\]
where $\widetilde{N}_{\delta}\left(r\right)=N_{\delta}\left(r\right)\cup\left\{ n\,\mid\, n<\sqrt{m\left(r\right)}\right\} $
and\[
\begin{array}{ll}
\mbox{{\rm \mbox{(ii)}}}_{n}: & |\xi_{n}|\le\frac{(a_{n}r^{n})^{-1}}{\sqrt{m\left(r\right)}},\\
\mbox{{\rm \mbox{(iii)}}}_{n}: & |\xi_{n}|\le\frac{1}{\sqrt{m\left(r\right)}},\\
\mbox{{\rm \mbox{(iv)}}}_{n}: & |\xi_{n}|\le e^{\frac{\delta n}{2}}.\end{array}\]
We notice that by \eqref{eq:m(r)_normal} and \eqref{eq:m(r)_low_bnd}
we have that $\#\widetilde{N}_{\delta}\left(r\right)\le C\sqrt{m\left(r\right)}$,
for $r$ large enough.
\begin{lem}
\label{lem:Low-Bnd-Est-For-Func}If $\Omega_{r}$ holds for a normal
value of $r$, then $f$ has no zeros inside $\disk r$.\end{lem}
\begin{proof}
We remind that $\delta=m\left(r\right)^{-1/4}$. To see that $f(z)$
has no zeros inside $\disk r$ we note that \begin{equation}
|f(z)|\ge|\xi_{0}|-\sum_{n=1}^{\infty}|\xi_{n}|a_{n}r^{n}.\label{eq:f_low_bnd}\end{equation}
First, we estimate the sum over the terms in $N(r)\backslash\left\{ 0\right\} $,
\[
\sum_{n\in N\left(r\right)\backslash\left\{ 0\right\} }|\xi_{n}|a_{n}r^{n}\le\sum_{n\in N\left(r\right)}\frac{1}{\sqrt{m\left(r\right)}}\le c_{1},\]
by \eqref{eq:m(r)_low_bnd}. Second, we estimate the sum over the
terms in $\widetilde{N}_{\delta}\left(r\right)\backslash N\left(r\right)$,
(notice that here $b_{n}\le0$ and $a_{n}r^{n}=e^{\log\left(a_{n}r^{n}\right)}=e^{nb_{n}}$)\[
\sum_{n\in\widetilde{N}_{\delta}\left(r\right)\backslash N\left(r\right)}|\xi_{n}|a_{n}r^{n}=\sum_{n\in\widetilde{N}_{\delta}\left(r\right)\backslash N\left(r\right)}|\xi_{n}|e^{nb_{n}}\le\sum_{n\in\widetilde{N}_{\delta}\left(r\right)\backslash N\left(r\right)}\frac{1}{\sqrt{m\left(r\right)}}\le c_{2},\]
using \eqref{eq:m(r)_normal}. Now the rest of the tail is bounded
by: (here $b_{n}\le-\delta$)\begin{eqnarray*}
\sum_{n\in\left(\widetilde{N}_{\delta}\left(r\right)\right)^{c}}|\xi_{n}|a_{n}r^{n} & = & \sum_{n\in\left(\widetilde{N}_{\delta}\left(r\right)\right)^{c}}|\xi_{n}|e^{nb_{n}}\le\sum_{n\in\left(\widetilde{N}_{\delta}\left(r\right)\right)^{c}}e^{\frac{\delta n}{2}}\cdot e^{-\delta n}\le\\
 & \le & \sum_{n\ge0}e^{-\frac{\delta n}{2}}\le\frac{3}{\delta}\le3\cdot\left(m\left(r\right)\right)^{1/4}.\end{eqnarray*}
So overall from \eqref{eq:f_low_bnd}\[
|f(z)|>C\left(m\left(r\right)\right)^{1/4}-c_{1}-c_{2}-3\left(m\left(r\right)\right)^{1/4},\]
if we select $C=4$ then for $r$ large enough we have that $f(z)\ne0$
inside $\disk r$. \end{proof}
\begin{lem}
\label{lem:Prob-of-Low-Bnd-Event}The probability of the event $\Omega_{r}$
is bounded from below by\textup{\[
\log\pr{\Omega_{r}}\ge-S(r)-C\cdot\sqrt{m\left(r\right)}\log m(r),\]
}for normal values of $r$ which are large enough.\end{lem}
\begin{proof}
By the properties of $\xi_{n}$ (see \ref{sub:Gauss-props}) we get\[
\pr{\mbox{{\rm \mbox{(i)}}}}=\exp\left(-C^{2}\cdot\sqrt{m\left(r\right)}\right)\]
and\[
\pr{\mbox{{\rm \mbox{(ii)}}}_{n}}\ge\frac{(a_{n}r^{n})^{-2}}{2m\left(r\right)},\]
therefore\begin{eqnarray*}
\pr{\mbox{(ii)}} & \ge & \prod_{n\in N\left(r\right)}\frac{(a_{n}r^{n})^{-2}}{2m\left(r\right)}=\prod_{n\in N\left(r\right)}e^{-2\cdot nb_{n}}\cdot\exp\left(-n\left(r\right)\log\left(2m\left(r\right)\right)\right)\ge\\
 & \ge & \exp\left(-2\cdot\sum_{n\in N\left(r\right)}nb_{n}\right)\cdot\exp\left(-Cn\left(r\right)\log m\left(r\right)\right)\ge\\
 & \ge & \exp\left(-S\left(r\right)-C\sqrt{m\left(r\right)}\log m\left(r\right)\right).\end{eqnarray*}
Similarly we have\[
\pr{\mbox{{\rm \mbox{(iii)}}}_{n}}\ge\frac{1}{2m\left(r\right)}\]
and so (by \eqref{eq:m(r)_normal} and \eqref{eq:m(r)_low_bnd})\[
\pr{\mbox{{\rm \mbox{(iii)}}}}\ge\exp\left(-C\sqrt{m\left(r\right)}\log m\left(r\right)\right).\]
\newpage{}Finally we have

\[
\pr{\mbox{{\rm \mbox{(iv)}}}_{n}^{c}}=\exp\left(-e^{\delta n}\right),\]
we use the following inequality (for some positive sequence $\left\{ A_{n}\right\} $)\[
\pr{\forall n\,:\,|\xi_{n}|\le A_{n}}=1-\pr{\exists n\,:\,|\xi_{n}|>A_{n}}\ge1-\sum\pr{|\xi_{n}|>A_{n}}.\]
Now\begin{eqnarray*}
\pr{\mbox{{\rm \mbox{(iv)}}}} & \ge & 1-\sum_{n\in\left(\widetilde{N}_{\delta}\left(r\right)\right)^{c}}\exp\left(-e^{\delta n}\right)\ge\\
 & \ge & 1-\sum_{n\ge\sqrt{m\left(r\right)}}\exp\left(-e^{\delta n}\right)\ge\\
 & \ge & 1-\sum_{n\ge\sqrt{m\left(r\right)}}\exp\left(-\delta n\right)\ge\\
 & \ge & 1-\frac{C}{\delta}\cdot\exp\left(-\frac{1}{\delta}\right)\ge\frac{1}{2},\end{eqnarray*}
for $r$ large enough (since $\sqrt{m\left(r\right)}=\frac{1}{\delta^{2}}$).
Since all the events are independent, in total the probablity is bounded
by:\begin{eqnarray*}
\pr{\Omega_{r}} & = & \pr{\mbox{{\rm \mbox{(i)}}}}\cdot\pr{\mbox{{\rm \mbox{(ii)}}}}\cdot\pr{\mbox{{\rm \mbox{(iii)}}}}\cdot\pr{\mbox{{\rm \mbox{(iv)}}}}\ge\\
 & \ge & \exp\left(-S\left(r\right)-C\sqrt{m\left(r\right)}\log m\left(r\right)\right).\end{eqnarray*}

\end{proof}
Proposition \ref{prp:p_H_upp_bnd} now follows from the previous lemmas.

\section{Bounds for Gauusian Entire Functions}

In this section we get some bounds for the modulus and logarithmic
derivative of Gaussian entire functions, which hold with high probability
(we use this term for events where the exceptional set is of small
probability in relation to the hole probability). These results will
be used in the next section in the proof of the lower bound.

\subsection{Bounds on the modulus of Gaussian entire functions}

We first bound the probability of the events when $M\left(r\right)$
is relatively large or small.
\begin{lem}
\textup{\label{lem:upp_bnd_M(r)}For normal values of $r$ which are
large enough, we have \[
\pr{M\left(r\right)\ge e^{3S\left(r\right)}}\le e^{-S^{2}\left(r\right)}.\]
}\end{lem}
\begin{proof}
We use the notation $\widetilde{N}\left(r\right)=N_{\delta}\left(r\right)\cup\left\{ n<S^{2}\left(r\right)\right\} $.
The proof is similar to the one of Proposition \ref{prp:p_H_upp_bnd}.
We define the following event $E$, which is the intersection of the
events $\mbox{{\rm \mbox{(i)}}}$ and $\mbox{(ii)}$\[
\begin{array}{ll}
\mbox{{\rm \mbox{(i)}}}: & {\displaystyle \bigcap_{n\in\tilde{N}\left(r\right)}}\mbox{{\rm \mbox{(i)}}}_{n},\\
\mbox{{\rm \mbox{(ii)}}}: & {\displaystyle \bigcap_{n\in\left(\widetilde{N}\left(r\right)\right)^{c}}}\mbox{{\rm \mbox{(ii)}}}_{n},\end{array}\]
and\[
\begin{array}{ll}
\mbox{{\rm \mbox{(i)}}}_{n}: & |\xi_{n}|\le\left(a_{n}r^{n}\right)^{-1}e^{2S\left(r\right)},\\
\mbox{{\rm \mbox{(ii)}}}_{n}: & |\xi_{n}|\le\exp\left(\frac{1}{2}\delta n\right).\end{array}\]
We have the following estimate for $M\left(r\right)$:\begin{eqnarray*}
|f(z)| & \le & \sum_{n\in\tilde{N}\left(r\right)}|\xi_{n}|a_{n}r^{n}+\sum_{n\in\left(\widetilde{N}\left(r\right)\right)^{c}}|\xi_{n}|a_{n}r^{n}\le\\
 & \le & \#\tilde{N}\left(r\right)\cdot e^{2S\left(r\right)}+\sum_{n\ge S^{2}\left(r\right)}e^{\frac{\delta n}{2}}\cdot e^{-\delta n}\le\\
 & \le & C\cdot S^{2}\left(r\right)\cdot e^{2S\left(r\right)}+\frac{C}{\delta}\cdot e^{-\frac{\delta S^{2}\left(r\right)}{2}}\le\\
 & \le & e^{3S\left(r\right)},\end{eqnarray*}
for $r$ large enough.

Now we estimate the probability of the complement of $E$. We have:\begin{eqnarray*}
\pr{|\xi_{n}|\ge\frac{e^{2S\left(r\right)}}{a_{n}r^{n}}} & = & \exp\left(-\frac{e^{4S\left(r\right)}}{\left(a_{n}r^{n}\right)^{2}}\right)\le\exp\left(-e^{2S\left(r\right)}\right),\\
\pr{|\xi_{n}|\ge e^{\frac{\delta n}{2}}} & = & \exp\left(-\exp\left(\delta n\right)\right).\end{eqnarray*}
By the union bound:\begin{eqnarray*}
\pr{\mbox{{\rm \mbox{(i)}}}^{c}} & \le & CS^{2}\left(r\right)\cdot\exp\left(-\exp\left(2S\left(r\right)\right)\right),\\
\pr{\mbox{{\rm \mbox{(ii)}}}^{c}} & \le & \sum_{n\ge S^{2}\left(r\right)}\exp\left(-\exp\left(\delta n\right)\right)\le\\
 & \le & C\cdot\exp\left(-\exp\left(S\left(r\right)\right)\right).\end{eqnarray*}
So overall\[
\pr{M\left(r\right)\ge e^{3S\left(r\right)}}\le\exp\left(-S^{2}\left(r\right)\right).\]

\end{proof}
In the other direction we have the following
\begin{lem}
\label{lem:low_bnd_M(r)}The probability of deviation from the lower
bound can be bounded by\[
\pr{M\left(r\right)\le\exp\left(-S\left(r\right)\right)}\le\exp\left(-S\left(r\right)\cdot n\left(r\right)\right).\]
\end{lem}
\begin{proof}
By Cauchy's estimate:\[
|\xi_{n}|a_{n}r^{n}\le M\left(r\right)\le e^{-S\left(r\right)}.\]
For $n\in N(r)$ we have\[
\pr{|\xi_{n}|\le\left(a_{n}r^{n}\right)^{-1}e^{-S\left(r\right)}}\le e^{-2S\left(r\right)},\]
and we get\[
\pr{M\left(r\right)\le\exp\left(-S\left(r\right)\right)}\le\prod_{n\in N(r)}e^{-2S\left(r\right)}\le\exp\left(-S(r)\cdot n\left(r\right)\right).\]

\end{proof}
Notice that there are no assumptions on (the normality of) $r$.

\subsection{Bounds for the logarithmic derivative}

In this section we assume that $f\left(z\right)\ne0$ inside $\disk R$,
and therefore $\log\left|f\right|$ is harmonic there. First we find
a bound for the average value of $\left|\log\left|f\right|\right|$.
$m$ denotes the normalized angular measure on $\cir r$. Under this
conditions we have
\begin{lem}
\label{lem:approx_log_int}For normal values of $R$ and outside an
exceptional set of probability at most\[
2\cdot\exp\left(-S(r)\cdot n\left(r\right)\right),\]
we have\[
\intop_{\cir R}\left|\log|f|\right|\, dm\le C\left(1-\frac{r}{R}\right)^{-2}\cdot S\left(R\right).\]
\end{lem}
\begin{proof}
Denote by $P_{j}(z)=P(z,z_{j})$ the Poisson kernel for the disk $\disk R,$
$|z|=R$, $|z_{j}|<R$. Using Lemma \ref{lem:low_bnd_M(r)}, we may
suppose that there is a point $a\in\cir r$ such that $\log|f(a)|\ge-S\left(r\right)$
(discarding an exceptional event of probability at most $\exp\left(-S\left(r\right)\cdot n\left(r\right)\right)$).
Then we have\[
-S\left(r\right)\le\intop_{\cir R}P(z,a)\log|f(z)|\, dm(z),\]
and hence\[
\intop_{\cir R}P(z,a)\log^{-}|f(z)|\, dm(z)\le\intop_{\cir R}P(z,a)\log^{+}|f(z)|\, dm(z)+S\left(r\right).\]
For $|z|=R$ and $|a|=r$ we have,\[
\frac{R-r}{2R}\le\frac{R-r}{R+r}\le P(z,a)\le\frac{R+r}{R-r}\le\frac{2R}{R-r}.\]
By Lemma \ref{lem:upp_bnd_M(r)}, outside a very small exception set
(of the order $e^{-S^{2}\left(R\right)}$), we have $\log M(R)\le3\cdot S(R)$.
Therefore\[
\intop_{\cir R}\log^{+}|f|\, d\mu\le3\cdot S(R).\]
Now we have\[
\intop_{\cir R}\log^{-}|f|\, d\mu\le\frac{2R}{R-r}\cdot S\left(r\right)+\frac{12R^{2}}{\left(R-r\right)^{2}}\cdot S\left(R\right).\]
Finally we get \begin{equation}
\intop_{\cir R}\left|\log|f|\right|\, d\mu\le\frac{CR^{2}}{\left(R-r\right)^{2}}\cdot S(R)\le C\left(1-\frac{r}{R}\right)^{-2}\cdot S(R)\label{eq:abs_log_int_upper_bnd}\end{equation}

\end{proof}
\newpage{}Now we find an upper bound for the (angular) logarithmic
derivative of $\log\left|f\right|$ inside $\disk R$. 
\begin{lem}
\label{lem:upp_bnd_log_deriv}Let $r<R$, then\[
\left|\frac{d\log\left|f\left(re^{i\phi}\right)\right|}{d\phi}\right|\le C\left(1-\frac{r}{R}\right)^{-5}\cdot S\left(R\right),\]
for normal values of $R$ and outside an exceptional set of probability
at most\[
2\cdot\exp\left(-S(r)\cdot n\left(r\right)\right).\]
\end{lem}
\begin{proof}
We start with the equation\[
\log\left|f\left(re^{i\phi}\right)\right|=\intop_{0}^{2\pi}\frac{R^{2}-r^{2}}{\left|Re^{i\theta}-re^{i\phi}\right|^{2}}\cdot\log\left|f\left(Re^{i\theta}\right)\right|\,\frac{d\theta}{2\pi},\]
so taking the derivative under the integral we get:\[
\frac{d\log\left|f\left(re^{i\phi}\right)\right|}{d\phi}=\intop_{0}^{2\pi}\frac{rR\left(R^{2}-r^{2}\right)\sin\left(\theta-\phi\right)}{\left|Re^{i\theta}-re^{i\phi}\right|^{4}}\cdot\log\left|f\left(Re^{i\theta}\right)\right|\,\frac{d\theta}{2\pi},\]
taking absolute value:\[
\left|\frac{d\log\left|f\left(re^{i\phi}\right)\right|}{d\phi}\right|\le\frac{C\left(R+r\right)}{\left(R-r\right)^{3}}\intop_{0}^{2\pi}\left|\log\left|f\left(Re^{i\theta}\right)\right|\right|\,\frac{d\theta}{2\pi},\]
using the previous lemma we get the required result.
\end{proof}

\section{Lower Bound for $p_{H}(r)$ \label{sec:Main-Thm-Lower-Bound}}

In order to find the lower bound for $p_{H}\left(r\right)$ we now
assume that $f(z)\ne0$ inside $\disk r$, for normal value of $r$.
We choose $\rho<r$, and write $\gamma=1-\frac{\rho}{r}$, where $\gamma$
will be small, depending on $r$.

The function $\log\left|f\left(z\right)\right|$ is harmonic in $\disk r$.
Therefore for $\rho<r$\[
\log\left|f\left(0\right)\right|=\intop_{0}^{2\pi}\log\left|f\left(\rho e^{i\alpha}\right)\right|\,\frac{d\alpha}{2\pi}.\]
Now, if we select $n$ points $z_{j}=\rho e^{i\theta_{j}}$, we have:\begin{eqnarray*}
\intop_{0}^{2\pi}\frac{1}{n}\sum_{j=1}^{n}\log\left|f\left(\rho e^{i\theta_{j}}\cdot e^{i\alpha}\right)\right|\,\frac{d\alpha}{2\pi} & = & \frac{1}{n}\sum_{j=1}^{n}\intop_{0}^{2\pi}\log\left|f\left(\rho e^{i\theta_{j}}\cdot e^{i\alpha}\right)\right|\,\frac{d\alpha}{2\pi}=\\
 & = & \frac{1}{n}\sum_{j=1}^{n}\log\left|f\left(0\right)\right|=\\
 & = & \log\left|f\left(0\right)\right|.\end{eqnarray*}
Since $\log\left|f\left(z\right)\right|$ is continuous inside $\disk r$,
we conclude that there exists some $\alpha^{*}$ such that\[
\frac{1}{n}\sum_{j=1}^{n}\log\left|f\left(\rho e^{i\theta_{j}}\cdot e^{i\alpha^{*}}\right)\right|=\log\left|f\left(0\right)\right|.\]
By Lemma \ref{lem:upp_bnd_log_deriv} the logarithmic derivative of
$f$ is not too large with high probability. Therefore, if $\alpha$
satifies \[
\left|\alpha-\alpha^{*}\right|<\Delta\alpha=\frac{c\gamma^{5}}{S\left(r\right)}\]
then:\[
\frac{1}{n}\sum_{j=1}^{n}\log\left|f\left(\rho e^{i\theta_{j}}\cdot e^{i\alpha}\right)\right|\le\log\left|f\left(0\right)\right|+1.\]
In this section we will show that the probability of the previous
event is very small. In particular we prove
\begin{prop}
\label{prp:p_H_low_bnd}For normal values of $r$, we have \[
p_{H}(r)\ge S(r)-Cn\left(r\right)\log S\left(r\right),\]
with $C$ some positive absolute constant.
\end{prop}

\subsection{Reduction to an estimate of a multivariate Gaussian event}

In this section we reduce the problem to an estimate of a probability
of an event with multivariate (complex) Gaussian distribution. We
first note that we work in the product space $\left\{ \left(\alpha,\omega\right)\in\left[0,2\pi\right]\times\Omega\right\} $,
where $\alpha$ is chosen uniformly in $\left[0,2\pi\right]$ and
$\Omega$ is the probability space for our Gaussian entire function
$f$ (we denote the probability measures by $m$ and $\mu$, respectively).
We define the following events (all depend on $r$):
\begin{itemize}
\item $H=\left\{ \left(\alpha,\omega\right)\,\mid\, f\left(z\right)\ne0\mbox{ in }\disk r\right\} $

\begin{itemize}
\item The hole event.
\end{itemize}
\item $L=\left\{ \left(\alpha,\omega\right)\,\mid\,\left|\frac{d\log\left|f\left(\rho e^{i\phi}\right)\right|}{d\phi}\right|\le C\cdot\left(1-\frac{\rho}{r}\right)^{-5}\cdot S\left(r\right),\,\forall\phi\in\left[0,2\pi\right]\right\} $ 

\begin{itemize}
\item The non-exceptional event of Lemma \ref{lem:upp_bnd_log_deriv}, w.r.t
$r$ and $\rho$
\end{itemize}
\item $C=\left\{ \left(\alpha,\omega\right)\,\mid\,\frac{1}{n}\sum\log\left|f\left(z_{j}e^{i\alpha}\right)\right|\le\log\left|f\left(0\right)\right|+1\right\} $ 
\item $D=\left\{ \left(\alpha,\omega\right)\,\mid\,\left|\alpha-\alpha^{\star}\left(\omega\right)\right|<\Delta\alpha\right\} $
\end{itemize}
We note that $\alpha^{\star}$ is measurable with respect to $\Omega$,
also the events $H$ and $L$ do not depend on the choice of $\alpha$.

By the previous discussion, we have:\[
\ind C{\alpha,\omega}\ge\ind{H\cap L\cap D}{\alpha,\omega}=\ind{H\cap L}{\alpha,\omega}\cdot\ind D{\alpha,\omega}=\ind{H\cap L}{\omega}\cdot\ind D{\alpha,\omega},\]
where $\ind A{\alpha,\omega}$ is the indicator function of the event
$A$. We now have by Fubini\begin{eqnarray*}
\pr C & =\\
\iint\ind C{\alpha,\omega}\, dm\left(\alpha\right)d\mu\left(\omega\right) & \ge & \int\left[\intop_{\left|\alpha-\alpha^{\star}\left(\omega\right)\right|<\Delta\alpha}1\, dm\left(\alpha\right)\right]\cdot\ind{H\cap L}{\omega}\, d\mu\left(\omega\right)=\\
 & = & 2\Delta\alpha\cdot\int\ind{H\cap L}{\omega}\, d\mu\left(\omega\right)\ge\Delta\alpha\cdot\pr{H\cap L}\ge\\
 & \ge & \Delta\alpha\cdot\left(\pr H-\pr{L^{c}}\right),\end{eqnarray*}
and so\[
\pr H\le\frac{1}{\Delta\alpha}\cdot\pr C+\pr{L^{c}}.\]
\newpage{}We now use the following events:
\begin{itemize}
\item $A=\left\{ \left(\alpha,\omega\right)\,\mid\,\log\left|f\left(0\right)\right|\le\log S\left(r\right)\right\} $
\item $B=\left\{ \left(\alpha,\omega\right)\,\mid\,{\displaystyle M\left(r\right)=\max_{z\in\disk r}\left|f\left(z\right)\right|}\le e^{3S\left(r\right)}\right\} $
\end{itemize}
To estimate $\pr C$ from above we use:\[
\pr C\le\pr{A\cap B\cap C}+\pr{A^{c}}+\pr{B^{c}}.\]
By \eqref{eq:Gaus_prob_large} and Lemma \ref{lem:upp_bnd_M(r)},
we have\begin{eqnarray*}
\pr{A^{c}} & \le & e^{-S^{2}\left(r\right)},\\
\pr{B^{c}} & \le & e^{-cS^{2}\left(r\right)},\end{eqnarray*}
note that these probabilities are very small with respect to the main
term.

In the next section, we will show that for a good selection of the
set of points $\left\{ z_{j}\right\} $ we have\begin{equation}
\pr{A\cap B\cap C}\le\exp\left(-S\left(\rho\right)+Cn\left(r\right)\log S\left(r\right)\right).\label{eq:low_bnd_main_est}\end{equation}
Therefore, we get (using Lemma \ref{lem:upp_bnd_log_deriv}),\begin{eqnarray*}
\pr H & \le & \pr{L^{c}}+\frac{1}{\Delta\alpha}\cdot\left(\pr{A\cap B\cap C}+\pr{A^{c}}+\pr{B^{c}}\right)\le\\
 & \le & 2e^{-S\left(\rho\right)\cdot n\left(\rho\right)}+\exp\left(-S\left(\rho\right)+\log\frac{1}{\Delta\alpha}+Cn\left(r\right)\log S\left(r\right)+\O\left(1\right)\right)\le\\
 & \le & \exp\left(-S\left(\rho\right)+c\log\frac{1}{\gamma}+Cn\left(r\right)\log S\left(r\right)+\O\left(1\right)\right).\end{eqnarray*}
Finally, by Lemma \ref{lem:S(r)_growth}, if we select $\gamma=\frac{1}{m\left(r\right)}$,
we have\[
\pr H\le\exp\left(-S\left(r\right)+Cn\left(r\right)\log S\left(r\right)\right),\]
thus proving Proposition \ref{prp:p_H_low_bnd}.

\subsection{Estimates for the probabilities}

We now turn to find an upper bound for the probability of the event
$A\cap B\cap C$, defined in the previous section, that is the event
when:\begin{eqnarray*}
\log\left|f\left(0\right)\right| & \le & \log S\left(r\right),\\
M\left(r\right)=\max_{z\in\disk r}\left|f\left(z\right)\right| & \le & e^{3S\left(r\right)},\\
\frac{1}{n}\sum\log\left|f\left(z_{j}e^{i\alpha}\right)\right| & \le & \log\left|f\left(0\right)\right|+1\le\log S\left(r\right)+1.\end{eqnarray*}
We remind that the points $\left\{ z_{j}\right\} $ were up until
now some $n$ arbitrary points on $\cir{\rho}$.

The vector $\left(f(z_{1}),\ldots,f(z_{n})\right)$ is distributed
like a multivariate complex Gaussian distribution, with the covariance
matrix:\[
\Sigma_{ij}=\Cov{f(z_{i}),f(z_{j})}=\Ex(f(z_{i})\overline{f(z_{j})})=\sum a_{k}^{2}\left(z_{i}\bar{z_{j}}\right)^{k}.\]
Note that the covariance matrix is invariant with respect to rotations
$f\left(e^{i\alpha}z\right)$. By Fubini we get\begin{eqnarray*}
\pr{A\cap B\cap C} & = & \iint\ind{A\cap B\cap C}{\alpha,\omega}\, dm\left(\alpha\right)d\mu\left(\omega\right)=\\
 & = & \intop\left[\intop\ind{A\cap B\cap C}{\alpha,\omega}\, d\mu\left(\omega\right)\right]\, dm\left(\alpha\right)\\
 & \le & \intop_{E}\frac{1}{\pi^{n}\det\Sigma}\exp\left(-\zeta^{*}\Sigma^{-1}\zeta\right)\, d\zeta,\end{eqnarray*}
where $E$ is the following set\textbackslash{}event:\[
E=\left\{ \left(f\left(z_{1}\right),\ldots,f\left(z_{n}\right)\right)\,\mid\,\frac{1}{n}\sum_{j=1}^{n}\log\left|f\left(z_{j}\right)\right|\le\log S\left(r\right)+1,\,\left|f\left(z_{j}\right)\right|\le e^{3S(r)},\,1\le j\le n\right\} .\]
We have the following upper bound for the probability of the event
$E:$\[
\pr E=\intop_{E}\frac{1}{\pi^{n}\det\Sigma}\exp\left(-\zeta^{*}\Sigma^{-1}\zeta\right)\, d\zeta\le\intop_{E}\frac{1}{\pi^{n}\det\Sigma}\, d\zeta=\frac{\volcn E}{\pi^{n}\det\Sigma}.\]
We start by finding a good lower bound for $\det\Sigma$, this depends
on a special selection of the points $\left\{ z_{j}\right\} $. 
\begin{lem}
Let $j_{1},\ldots,j_{n-1}\in\bbn$. There exist $n$ points $\left\{ z_{j}\right\} $
on $\cir r$ such that the determinant of the following generalized
Vandermonde matrix \[
A=\left(\begin{array}{cccc}
1 & z_{1}^{j_{1}} & \ldots & z_{1}^{j_{n-1}}\\
\vdots & \vdots & \vdots & \vdots\\
1 & z_{n}^{j_{1}} & \ldots & z_{n}^{j_{n-1}}\end{array}\right)\]
satisfies $\left|\det A\right|\ge r^{\sum_{k=1}^{n-1}j_{k}}$.\end{lem}
\begin{proof}
Start with the formal expression for the determinant\[
\det A=\sum_{\sigma}\mathrm{sgn}\left(\sigma\right)\prod_{m=0}^{n-1}z_{m}^{j_{\sigma\left(m\right)}}.\]
Write $z_{j}=re^{i\theta_{j}}$, now we have:\begin{eqnarray*}
\intop_{\mathbb{T}^{n}}\left|\det A\right|^{2}\, d\theta_{1}\ldots d\theta_{n} & = & \intop_{\mathbb{T}^{n}}\sum_{\sigma}\mathrm{sgn}\left(\sigma\right)\prod_{m=1}^{n}z_{m}^{j_{\sigma\left(m\right)}}\cdot\sum_{\sigma}\mathrm{sgn}\left(\sigma\right)\prod_{m=1}^{n}\bar{z}_{m}^{j_{\sigma\left(m\right)}}\, dz_{1}\ldots dz_{n}=\\
 & = & r^{2\cdot\sum_{k=1}^{n-1}j_{k}}\cdot\left(\left(2\pi\right)^{n}\cdot\sum_{\sigma}1+\intop_{\mathbb{T}^{n}}\sum_{\sigma\ne\tau}\mathrm{sgn}\left(\sigma\right)\cdot\mathrm{sgn}\left(\tau\right)e^{i\cdot\sum a_{\sigma,\tau,j}\theta_{j}}\right)\ge\\
 & \ge & r^{2\cdot\sum_{k=1}^{n-1}j_{k}}\cdot n!,\end{eqnarray*}
since at least one of the $a_{\sigma,\tau,j}$ is different than $0$
in each sum.
\end{proof}
We now have the following
\begin{cor}
\label{cor:low_bnd_det_sigma}Using the configuration of the points
$\left\{ z_{j}\right\} $ given in the previous lemma, we have\[
\log\left(\det\Sigma\right)\ge S(r).\]
\end{cor}
\begin{proof}
Notice that we can represent $\Sigma$ in the following form\[
\Sigma=V\cdot V^{*}\]
where \[
V=\left(\begin{matrix}a_{0} & a_{1}\cdot z_{1} & \dots & a_{n}\cdot z_{1}^{n} & \dots\\
\vdots & \vdots & \vdots & \vdots & \dots\\
a_{0} & a_{1}\cdot z_{N} & \dots & a_{n}\cdot z_{N}^{n} & \dots\end{matrix}\right).\]
We can estimate the determinant of $\Sigma$ by projecting $V$ on
$N(r)=\left\{ a_{0},a_{j_{1}},\ldots,a_{j_{n-1}}\right\} $ coordinates
(let's denote this projection by $P$). Since $\det\Sigma$ is the
square of the product of the singular values of $V$, and these values
are only reduced by the projection, we have\[
\det\Sigma\ge\left(\det PV\right)^{2}=\left|\begin{array}{cccc}
a_{0} & a_{j_{1}}z_{1}^{j_{1}} & \ldots & a_{j_{n-1}}z_{1}^{j_{n-1}}\\
\vdots & \vdots & \vdots & \vdots\\
a_{0} & a_{j_{1}}z_{n}^{j_{1}} & \ldots & a_{j_{n-1}}z_{n}^{j_{n-1}}\end{array}\right|^{2}\]
and so by the previous lemma\begin{eqnarray*}
\det\Sigma & \ge & \prod_{j\in N_{1}\left(r\right)}a_{j}^{2}\cdot\left|\begin{array}{cccc}
1 & z_{1}^{j_{1}} & \ldots & z_{1}^{j_{n-1}}\\
\vdots & \vdots & \vdots & \vdots\\
1 & z_{n}^{j_{1}} & \ldots & z_{n}^{j_{n-1}}\end{array}\right|^{2}\\
 & \ge & \prod_{j\in N_{1}\left(r\right)}a_{j}^{2}\cdot r^{2\cdot\sum_{j\in N_{1}\left(r\right)}j}\\
 & = & \prod_{j\in N_{1}\left(r\right)}a_{j}^{2}r^{2j}=\exp\left(S\left(r\right)\right)\end{eqnarray*}

\end{proof}
We now want to estimate the integral\[
I=\intop_{E}\frac{1}{\pi^{n}}\, d\zeta,\]
that is to estimate the volume of the following set:\[
E=\left\{ \zeta\in\bbc^{n}\,\mid\,\frac{1}{n}\sum_{j=1}^{n}\log\left|\zeta_{j}\right|\le\log S\left(r\right)+1\mbox{ and}\left|\zeta_{j}\right|\le e^{3S\left(r\right)},\,1\le j\le n\right\} ,\]
with respect to the Lebesgue measure on $\bbc^{n}.$ We will use the
following lemma (see \cite[Lemma 11]{Ni1}) 
\begin{lem}
Set $s>0$, $t>0$ and $N\in\bbn^{+}$, such that $\log\left(t^{N}/s\right)\ge N$.
Denote by $\C_{N}$ the following set\[
\C_{N}=\C_{N}\left(t,s\right)=\left\{ r=\left(r_{1},\ldots,r_{N}\right)\::\:0\le r_{j}\le t,\,\prod_{1}^{N}r_{j}\le s\right\} .\]
Then\[
\volrn N{\C_{N}}\le\frac{s}{\left(N-1\right)!}\log^{N}\left(t^{N}/s\right).\]

\end{lem}
\newpage{}Now we have the almost immediate
\begin{cor}
\label{cor:prob_integral_upper_bnd}Suppose that $r$ is normal and
large enough, then we have\[
I\le\exp\left(Cn\left(r\right)\log S\left(r\right)\right).\]
\end{cor}
\begin{proof}
To shorten the expressions, we write\begin{eqnarray*}
n & = & n(r)\\
s & = & \exp\left(n\left(r\right)\left(\log S\left(r\right)+1\right)\right)\\
t & = & \exp\left(3S(r)\right).\end{eqnarray*}
We want to translate the integral $I$ into an integral in $\bbr^{N}$,
using the change of variables $\zeta_{j}=r_{j}\cos(\theta_{j})+ir_{j}\sin(\theta_{j})$.
Integrating out the variables $\theta_{j}$, we get $I^{\prime}=2^{n}\intop_{\C}\prod r_{j}\, dr$,
where the new domain is\[
\C=\left\{ r=\left(r_{1},\ldots,r_{n}\right)\::\:0\le r_{j}\le t,\,\prod_{j=1}^{n}r_{j}\le s\right\} .\]
We can find an explicit expression for this integral, but, instead
we will simplify it even more to\begin{equation}
I^{\prime}\le2^{n}s\cdot\volrn n{\C}\label{eq:I_prime_upp_bnd}\end{equation}
Now, in order to use the previous lemma, we have to check the condition
$\log\left(t^{n}/s\right)\ge n$, or (where $C>0$)\[
3n(r)S(r)-n(r)\left(\log S(r)+1\right)\ge n(r),\]
which is satisfied under our assumptions, for $r$ large enough. After
applying the lemma, we get (for $r$ large enough)\begin{eqnarray*}
I^{\prime} & \le & \frac{n\cdot2^{n}s^{2}}{n!}\log^{n}\left(t^{n}/s\right)\\
 & \le & \frac{s^{2}e^{2n}}{n^{n}}\log^{n}\left(t^{n}/s\right)\\
 & = & \exp\left(2\log s+n\log_{2}t+2n-n\log_{2}s\right)\\
 & \le & \exp\left(2\log s+n\log_{2}t\right).\end{eqnarray*}
Recalling the definitions of $n,s$ and $t$, we finally get

\begin{eqnarray*}
\log I^{\prime} & \le & 2n\left(r\right)\left(\log S(r)+1\right)+n\left(r\right)\log S\left(r\right)+Cn\left(r\right)\\
 & \le & Cn\left(r\right)\log S\left(r\right).\end{eqnarray*}

\end{proof}
Now, the estimate \eqref{eq:low_bnd_main_est} follows, since the
original points $\left\{ z_{j}\right\} $ satisfy $\left|z_{j}\right|=\rho$
and by Corollaries \ref{cor:low_bnd_det_sigma} and \ref{cor:prob_integral_upper_bnd}:\[
\pr E\le\exp\left(-S\left(\rho\right)+Cn\left(r\right)\log S\left(r\right)\right).\]
This completes the proof of Proposition \ref{prp:p_H_low_bnd}.


\begin{thebibliography}{ST}
\bibitem[BKPV]{BKPV}J. Ben Hough, M. Krishnapur, Y. Peres, B. Virág,
{}``Zeros of Gaussian Analytic Functions and Determinantal Point
Processes'', American Mathematical Society, 2010 (University lecture
series ; v. \textbf{51}).

\bibitem[Hay]{Hay}W. K. Hayman, The local growth of power series:
a survey of the Wiman-Valiron method. Canad. Math. Bull. \textbf{17}
(1974), no. 3, 317--358. 

\bibitem[Kah]{Kah}J.-P. Kahane, {}``Some random series of functions'',
2nd ed., Cambridge Studies in Advanced Mathematics, vol. \textbf{5},
1985.

\bibitem[Ni1]{Ni1}A. Nishry, Asymptotics of the hole probability
for zeros of random entire functions. Int. Math. Res. Not. IMRN 2010,
no. \textbf{15}, 2925\textendash{}2946. 

\bibitem[Ni2]{Ni2}A. Nishry, The Hole Probability for Gaussian Entire
Functions, to be published in Israel Journal of Mathematics, arXiv:
math.CV/0909.1270.

\bibitem[ST]{ST}M. Sodin and B. Tsirelson, Random Complex zeroes.
III. Decay of the hole probability. Israel J. Math \textbf{147} (2005),
371\textendash{}379.
\end{thebibliography}
\end{document}